# SURVIVAL AND COEXISTENCE FOR A MULTITYPE CONTACT PROCESS

By J. Theodore Cox[1] and Rinaldo B. Schinazi[2]

*Syracuse University and University of Colorado, Colorado Springs*

We study the ergodic theory of a multitype contact process with equal death rates and unequal birth rates on the $d$-dimensional integer lattice and regular trees. We prove that for birth rates in a certain interval there is coexistence on the tree, which by a result of Neuhauser is not possible on the lattice. We also prove a complete convergence result when the larger birth rate falls outside of this interval.

**1. Introduction and results.** We consider the multitype contact process on a countable set $S$, introduced by Neuhauser in [8]. In this paper $S$ is either the $d$-dimensional integer lattice $\mathbb{Z}^d$, the case considered by Neuhauser, or the homogeneous connected tree $\mathbb{T}_d$ in which each vertex has $d+1$ neighbors, $d \geq 2$. The primary reason we are interested in the multitype contact process defined on the tree $\mathbb{T}_d$ is that, as we show below, it exhibits phenomena that the process defined on $\mathbb{Z}^d$ does not.

In this model each point, or site, of $S$ is either vacant or occupied by an individual of one of two possible types. The system is described by a configuration $\xi \in \{0,1,2\}^S$, where $\xi(x) = 0$ means that site $x$ is vacant, and for $i = 1, 2$, $\xi(x) = i$ means that $x$ is occupied by an individual of type $i$. For $x, y \in S$ write $x \sim y$ if $x, y$ are nearest neighbors, and define

$$^i\xi = \{x : \xi(x) = i\}, \qquad \xi \in \{0,1,2\}^S, i = 1, 2.$$

For $x \in S$ and $\xi \in \{0,1,2\}^S$, let $n_i(x, \xi)$ denote the number of neighbors of $x$ that are of type $i$,

$$n_i(x, \xi) = \sum_{y \sim x} 1\{\xi(y) = i\}, \qquad i = 1, 2.$$

Received June 2007; revised April 2008.
[1]Supported in part by NSF Grant 0505439.
[2]Supported in part by NSF Grant 0701396.

The multitype contact process $\xi_t$ with birth rates $\lambda_1, \lambda_2$ is the Feller process taking values in $\{0,1,2\}^S$ which makes transitions at $x$ in configuration $\xi$

$$i \to 0 \qquad \text{at rate } 1, i = 1, 2, \tag{1.1}$$

$$0 \to i \qquad \text{at rate } \lambda_i n_i(x, \xi), i = 1, 2. \tag{1.2}$$

Informally, an individual of either type dies at rate one, and for each site $x$ and each neighboring site $y$, an individual of type $i$ at $x$ gives birth at rate $\lambda_i$ to an individual of the same type at $y$, provided that site is vacant. Thus, the two types interact in their competition for space. Existence and uniqueness of a Feller process $\xi_t$ determined by the above rates follows from general results of [6], Theorem I.3.9. We will give a "graphical construction" of $\xi_t$ in Section 2 below.

For this process on $\mathbb{Z}^d$, Neuhauser proved in [8] that the two types can coexist if and only if $d \geq 3$ and $\lambda_1 = \lambda_2 > \lambda_c$ ($\lambda_c$ is defined below). In [8] coexistence meant the existence of a translation invariant measure $\nu$ which is invariant for the process and concentrates on configurations with infinitely many individuals of each type. That is, letting $|\cdot|$ denote cardinality, $\nu$ should satisfy

$$\nu(\xi : |^1\xi| = |^2\xi| = \infty) = 1. \tag{1.3}$$

The noncoexistence part of this result was strengthened by Durrett and Neuhauser in [2]. By Theorem 2 there, if $\lambda_2 > \lambda_1 \vee \lambda_c$, and if the initial configuration $\xi_0$ has infinitely many type 2 individuals, then for every site $x$, $P(\xi_t(x) = 1) \to 0$ as $t \to \infty$. Consequently, there is no invariant measure $\nu$, translation invariant or otherwise, satisfying (1.3).

Switching from $\mathbb{Z}^d$ to $\mathbb{T}_d$, one can use the approach from [8] to show that if $\lambda_1 \neq \lambda_2$, then there can be no invariant measure $\nu$ which is homogeneous (invariant under the obvious shifts) and also satisfies (1.3). However, the arguments in [2], which involve a block construction, do not seem to be directly applicable on the tree, leaving open the possibility of nonhomogeneous invariant measures satisfying (1.3). We show in Theorem 1 that such measures do in fact exist, provided the birth rates lie in a certain interval. In Theorem 2 we show that there is no coexistence when the birth rates are outside this interval. For this result our method of proof applies equally well on $\mathbb{Z}^d$.

In order to state our results, it is necessary to briefly review the basic, single-type contact process $\zeta_t$ introduced by Harris in [4], treated in detail in Chapter VI of [6] and in Part I of [7]. In this process each site of $S$ is either vacant or occupied by an individual, and $\zeta_t$ is the set of occupied sites at time $t$. Supposing $S$ is either $\mathbb{Z}^d$ or $\mathbb{T}_d$ and $A \subset S$, the transition rates are given by

$$A \to A \setminus \{x\} \qquad \text{for } x \in A \text{ at rate } 1,$$

$$A \to A \cup \{x\} \qquad \text{for } x \notin A \text{ at rate } \lambda|\{y \in A : y \sim x\}|.$$

For $A \subset S$ let $\zeta_t^A$ denote the process with initial state $\zeta_0^A = A$, and for $x \in S$ write $\zeta_t^x$ for $\zeta_t^{\{x\}}$. We give a construction of $\zeta_t$ in Section 2 below.

Define the two critical values $\lambda_* \le \lambda^*$ by

$$\lambda_* = \inf\{\lambda : P(\zeta_t^x \ne \varnothing \ \forall t > 0) > 0\}$$

and

$$\lambda^* = \inf\{\lambda : P(x \in \zeta_t^x \text{ i.o. as } t \to \infty) > 0\}.$$

By translation invariance, the above probabilities do not depend on $x$. The contact process is said to *die out* if $P(|\zeta_t^x| \ge 1 \ \forall t > 0) = 0$, *survive weakly* if this probability is positive but $P(x \in \zeta_t^x \text{ i.o. as } t \to \infty) = 0$, and *survive strongly* if this last probability is positive. It is known that for $S = \mathbb{Z}^d$ weak survival does not occur, $\lambda_* = \lambda^*$, and the common value $\lambda_c$ is finite and strictly positive. However, for the tree $\mathbb{T}_d$, weak survival does occur. It is known that $0 < \lambda_* < \lambda^* < \infty$, the process dies out if $\lambda = \lambda_*$ and survives weakly if $\lambda = \lambda^*$ (see Proposition I.4.39 and Theorem I.4.65 in [7]). It turns out that these different possibilities for survival must be taken into account when investigating the question of coexistence for the multitype process on $\mathbb{T}_d$.

For the tree $\mathbb{T}_d$ fix a site $\mathcal{O} \in \mathbb{T}_d$ and call it the root. For $x \ne y \in \mathbb{T}_d$ there is a unique sequence of distinct sites $x = x_0, x_1, \ldots, x_n = y$ such that $x_{i-1} \sim x_i$ for $i = 1, \ldots, n$. Let this $n$ be the distance from $x$ to $y$, $|x - y| = n$. Let $B_K$ be the ball of radius $K$ centered at $\mathcal{O}$, $B_K = \{y \in \mathbb{T}_d : |y - \mathcal{O}| \le K\}$, and let $\partial B_K$ denote the outer boundary, $\partial B_K = \{y \in \mathbb{T}_d : |y - \mathcal{O}| = K + 1\}$. For $x \in \mathbb{T}_d$ let $S(x)$ be the sector of the tree starting at $x$ pointing away from the root, that is,

$$S(x) = \{y \in \mathbb{T}_d : |y - \mathcal{O}| = |y - x| + |x - \mathcal{O}|\}.$$

A nice set of configurations to work with is $\Xi_0$, the set of configurations $\eta \in \{0, 1, 2\}^{\mathbb{T}_d}$ which satisfy:

$$\exists K < \infty \qquad \text{such that } \eta \text{ is constant on each } S(y), y \in \partial B_K. \tag{1.4}$$

Our first results are for the case that both birth rates lie between the two contact process critical values. For this case there is coexistence, and even weak convergence starting from any configuration in $\Xi_0$. For a probability measure $\nu$ on $\{0, 1, 2\}^S$, $\xi_t \Rightarrow \nu$ as $t \to \infty$ means that the finite dimensional distributions of $\xi_t$ converge to those given by $\nu$.

THEOREM 1. *Assume $S = \mathbb{T}_d$ and $\lambda_1, \lambda_2 \in (\lambda_*, \lambda^*]$.*

(i) *If $\xi_0 \equiv i$ on $S(y)$ for some site $y$, then for all sites $x$,*

$$\liminf_{t \to \infty} P(\xi_t(x) = i) > 0. \tag{1.5}$$

(ii) *If* $\xi_0 = \eta \in \Xi_0$, *then*

$$\xi_t \Rightarrow \nu_\eta \qquad \text{as } t \to \infty \tag{1.6}$$

*for some measure* $\nu_\eta$ *which is necessarily invariant for* $\xi_t$.

(iii) *Assume* $\eta \in \Xi_0$. *If* $|{}^i\eta| = \infty$, *then* $\nu_\eta(\xi : |{}^i\xi| = \infty) = 1$. *If* $\eta \equiv i$ *on some sector* $S(y)$ *and* $|x| \to \infty$, $x \in S(y)$, *then* $\nu_\eta(\xi : \xi(x) = j) \to 0$ *if* $j$ *is not 0 or* $i$.

(iv) *Assume* $\eta, \eta' \in \Xi_0$. *If* $|\{x : \eta(x) \neq \eta'(x)\}| < \infty$, *then* $\nu_\eta = \nu_{\eta'}$. *If* $\exists y \in \mathbb{T}_d$ *and* $i \neq j$ *such that* $\eta \equiv i$ *and* $\eta' \equiv j$ *on* $S(y)$, *then* $\nu_\eta \neq \nu_{\eta'}$.

Theorem 1 shows there is coexistence even for unequal birth rates, provided the two rates lie in the contact process weak survival interval, and exhibits a large class of nonhomogeneous invariant measures $\nu_\eta$. If $\eta \in \Xi_0$ is identically 1 on some sector and identically 2 on another, then by (ii) and (iii) above, $\nu_\eta$ is a nonhomogeneous invariant measure which satisfies (1.3). This means that when $\lambda_1 < \lambda_2$ and both rates lie in $(\lambda_*, \lambda^*]$, the 2's are not strong enough to drive the 1's from bounded regions of the tree, and coexistence is possible. Two-type competition models have been studied by many others; see [1, 3] and [5] for instance. The model in [5] exhibits "global" coexistence with unequal rates, in that $P(|{}^1\xi_t| \geq 1, |{}^2\xi_t| \geq 1 \ \forall t > 0) > 0$ for every initial configuration with $|{}^1\xi_0| \geq 1, |{}^2\xi_0| \geq 1$. Coexistence results are proved in [1] for the multitype contact process on $\mathbb{Z}^d$ with long-range interactions, or an additional "death" mechanism. Theorem 1 may give the first "local" coexistence result for a nearest-neighbor interaction model with equal death rates and unequal birth rates.

It is a different story when $\lambda_1 \neq \lambda_2$ with one or both rates larger than $\lambda^*$. In this case coexistence, even in the weak sense of (1.7) below, is not possible.

THEOREM 2. *Assume* $S$ *is* $\mathbb{Z}^d$ *or* $\mathbb{T}_d$ *and* $\lambda_2 > \lambda_1 \vee \lambda^*$. *For any initial configuration* $\xi_0$ *and* $x \in S$,

$$\lim_{t \to \infty} P(\xi_t(x) = 1 \ \text{and} \ |{}^2\xi_t| \geq 1) = 0. \tag{1.7}$$

*Consequently,* $\lim_{t\to\infty} P(\xi_t(x) = 1) = 0$ *if* $|{}^2\xi_0| = \infty$.

For $S = \mathbb{Z}^d$, (1.7) follows from the results in Section 3 of [2]. We include this case in the statement of Theorem 2 since our proof for $S = \mathbb{T}_d$ applies equally well to the lattice case, and is simpler than the one in [2]. The lack of coexistence in (1.7) means that an invariant measures can concentrate on only one type, and allows us to prove that the process converges weakly from any initial configuration. We need some additional notation and information about the basic contact process before stating our results.

For the single-type contact process $\zeta_t$ with birth rate $\lambda$, define the survival probabilities $\alpha_A = \alpha_A(\lambda)$ by

$$\alpha_A = P(\zeta_t^A \neq \varnothing \ \forall t \geq 0), \qquad A \subset S, \tag{1.8}$$

and let $\alpha = \alpha_{\{x\}}$, which by translation invariance does not depend on $x$.

It is well known (see (I.1.4) of [7]) that there is a probability measure $\bar{\nu} = \bar{\nu}_\lambda$ called the upper invariant measure such that

$$\zeta_t^S \Rightarrow \bar{\nu}_\lambda \qquad \text{as } t \to \infty. \tag{1.9}$$

Letting $\delta_\varnothing$ be the unit point mass on the empty set, $\bar{\nu} \neq \delta_\varnothing$ if and only if $\lambda > \lambda_*$. The *complete convergence theorem* for $\zeta_t$ is the statement

$$\zeta_t^A \Rightarrow (1-\alpha_A)\delta_\varnothing + \alpha_A \bar{\nu} \qquad \text{as } t \to \infty \ \forall A \subset S. \tag{1.10}$$

It is known that (1.10) holds for both $\mathbb{Z}^d$ and $\mathbb{T}_d$ if $\lambda > \lambda^*$ (see Theorems I.2.27 and I.4.70 in [7]).

We show here that there is an analogous theorem for the multitype contact process if $\lambda_2 > \lambda_1 > \lambda^*$. To state it, we must define appropriate survival probabilities and "upper invariant measures." For birth rates $\lambda_1, \lambda_2$ and configurations $\eta \in \{0,1,2\}^S$, let $\xi_0 = \eta$ and define $\alpha_\eta^i = \alpha_\eta^i(\lambda_1, \lambda_2)$, $i = 1, 2$, by

$$\begin{aligned} \alpha_\eta^1 &= P({}^1\xi_t \neq \varnothing \ \forall t \geq 0 \text{ and } {}^2\xi_t = \varnothing \ \text{ eventually}), \\ \alpha_\eta^2 &= P({}^2\xi_t \neq \varnothing \ \forall t \geq 0). \end{aligned} \tag{1.11}$$

We need probability measures $\bar{\nu}^i = \bar{\nu}^i_{\lambda_i}$ on $\{0,1,2\}^S$, which correspond to $\bar{\nu}_{\lambda_i}$ and concentrate on configurations in which all individuals are of type $i$. These measures are defined by the requirements

$$\bar{\nu}^1(\xi : {}^2\xi = \varnothing) = \bar{\nu}^2(\xi : {}^1\xi(x) = \varnothing) = 1$$

and

$$\bar{\nu}^i(\xi : {}^i\xi \cap A \neq \varnothing) = \bar{\nu}_{\lambda_i}(\zeta : \zeta \cap A \neq \varnothing), \qquad A \subset S, i = 1, 2.$$

With these definitions in place, we can now state our complete convergence theorem for $\xi_t$ when $\lambda_2 > \lambda_1 > \lambda^*$. In the following let $\delta_{\mathbf{0}}$ denote the measure on $\{0,1,2\}^S$ which concentrates on the single configuration $\xi \equiv 0$.

THEOREM 3. *Assume $S$ is $\mathbb{Z}^d$ or $\mathbb{T}^d$, $\lambda_2 > \lambda_1 > \lambda^*$ and $\xi_0 = \eta$. Then*

$$\xi_t \Rightarrow (1 - \alpha_\eta^1 - \alpha_\eta^2)\delta_{\mathbf{0}} + \alpha_\eta^1 \bar{\nu}^1 + \alpha_\eta^2 \bar{\nu}^2 \qquad \textit{as } t \to \infty. \tag{1.12}$$

Given Theorem 2, this result is not surprising. If the 2's survive, then the 1's are driven out of bounded regions, so the 2's in effect form a single type contact process and (1.10) takes over. If there are finitely many 2's which die out while the 1's manage to survive, then the 1's form a single type contact process and again (1.10) takes over.

The complete convergence theorem for the contact process (1.10) does not hold on $\mathbb{T}_d$ for $\lambda \in (\lambda_*, \lambda^*]$. In this case the contact process has a wide variety of invariant measures (see Theorems I.4.107 and I.4.121 of [7]) and, hence, possible limits for $\zeta_t$. We cannot expect (1.12) to hold as stated if $\lambda_1 \le \lambda^* < \lambda_2$. However, if we restrict $\xi_0$ to configurations for which the corresponding contact process of 1's converges weakly, then $\xi_t$ also converges weakly.

THEOREM 4. *Assume $S = \mathbb{T}^d$, $\lambda_2 > \lambda^* \ge \lambda_1$ and $\xi_0 = \eta$. Let $\zeta_t$ be the single-type contact process with birth rate $\lambda_1$ and initial state $\zeta_0 = {}^1\eta$, and assume that $\zeta_t \Rightarrow \mu$ as $t \to \infty$. Then*

$$\xi_t \Rightarrow (1 - \alpha_\eta^2)\bar{\mu}^1 + \alpha_\eta^2 \bar{\nu}^2 \qquad \text{as } t \to \infty, \tag{1.13}$$

*where $\mu^1(\xi : {}^2\xi = \varnothing) = 1$ and $\mu^1(\xi : {}^1\xi \cap A \ne \varnothing) = \mu(\zeta : \zeta \cap A \ne \varnothing)$.*

For $S = \mathbb{Z}^d$, the conclusions of Theorems 2 and 3 can be derived from results in [2] (Lemma 3 and the construction in Section 3 there). However, our methods are simpler and should apply without much change to other choices for $S$, including some periodic graphs.

As previously noted, for $S = \mathbb{Z}^d$ and $\lambda_1 = \lambda_2 > \lambda_c$, it was shown in [8] that there is coexistence for $d \ge 3$ but not for $d \le 2$. Presumably the arguments for these results can be adapted to handle the tree $\mathbb{T}_d$, where one expects there should be coexistence for $\lambda_1 = \lambda_2 > \lambda^*$.

In the next section we construct our processes using the standard "graphical construction" via Poisson processes. The construction naturally contains various couplings and dual processes used in our proofs. In Sections 3–6 we prove Theorems 1–4, respectively.

We note that our main tool is the *ancestor duality* introduced by Neuhauser in [8].

**2. Construction and duality.** We start by constructing our process using Harris' graphical method, assuming from now on that

$$\lambda_1 \le \lambda_2. \tag{2.1}$$

The construction takes place in the space–time set $S \times [0, \infty)$ using independent families of Poisson processes. For $x \in S$ let $\{T_n^x : n \ge 1\}$ be the arrival times of a Poisson process with rate 1. At the times $T_n^x$ we put a $\delta$ at site $x$ to indicate that there is a death at $x$: if site $x$ is occupied by either type, it becomes vacant at that time. For all pairs of nearest neighbors $x,y \in S$ let $\{B_n^{x,y} : n \geq 1\}$ be the arrival times of a Poisson process with rate $\lambda_2$. At the times $B_n^{x,y}$ we do two things. We draw an arrow from site $x$ to site $y$, and with probability $1 - \lambda_1/\lambda_2$, independently of everything else, label the arrow with a "2" (and otherwise do not label the arrow). If there is a 2 at $x$ and $y$ is vacant at that time, then there is a birth of a 2 at $y$. If there is a 1 at $x$ and $y$ is vacant, we put a 1 at $y$ provided the arrow does not have a 2 on it. Thus, the 2-arrows are really "2-only" arrows. If $\lambda_1 = \lambda_2$, then no arrow is marked with a 2. The Poisson processes $T^x, B^{x,y}$, $x,y \in S$, are all independent of one another.

For sites $x,y \in S$ and times $0 < s \leq t$, we say there is a path up from $(x,s)$ to $(y,t)$ if there is a sequence of times $t_0 = s < t_1 < t_2 < \cdots < t_n = t$ and a sequence of sites $x_0 = x, x_1, \ldots, x_n = y$ such that, for $i = 1, 2\ldots, n$, $x_{i-1} \sim x_i$, there is an arrow from $x_{i-1}$ to $x_i$ at time $t_i$ and the time segments $\{x_i\} \times (t_{i-1}, t_i)$ do not contain any $\delta$'s. By default, there is always a path up from $(x,t)$ to $(x,t)$. A path up which has at least one arrow labeled 2 will be called a "2-path," and a path with no arrows labeled 2 will be called a "1-path." Note that 1's propagate *only* along 1-paths, but 2's propagate along both 1-paths and 2-paths. For $s < t$, there is an $i$-path down from $(y,t)$ to $(x,s)$ if and only if there is an $i$-path up from $(x,s)$ to $(y,t)$. Given an initial configuration $\xi_0$, we may construct $\xi_t, t \geq 0$, from our Poisson processes by following paths from occupied sites forward in time, using $\delta$'s for deaths and arrows for births, as appropriate.

We now define several "reverse" time dual processes, starting with the simplest. For $x \in S$, $t > 0$, $i = 1, 2$ and $0 \leq s \leq t$, define

$$D_s^{i,(x,t)} = \{y \in S : \text{ there is an } i\text{-path down from } (x,t) \text{ to } (y, t-s)\}$$

and $D_s^{(x,t)} = D_s^{1,(x,t)} \cup D_s^{2,(x,t)}$. The *ancestor process* introduced in [8] is more complicated. Fix $x \in S$ and $t > 0$, and consider $D_s^{(x,t)}, 0 \leq s \leq t$. If $D_s^{(x,t)}$ is not empty, then the sites in $D_s^{(x,t)}$ are the possible *ancestors* at (forward) time $t - s$ of $(x,t)$, which can be arranged in decreasing order of priority, $(a_1(s), a_2(s), \ldots, a_n(s))$ for some $n$, with $a_1(s)$ denoting the *primary* ancestor (see [8]). The $j$th ancestor $a_j$ is associated with a path up from $(a_j, t-s)$ to $(x,t)$, which may or may not contain an arrow labeled 2, blocking propagation of 1's. We will use an equivalent but slightly different formulation of this process, which we now describe.

An ancestor configuration $\widehat{\xi}$ is either the empty set, or a sequence of pairs $((a_1, b_1), \ldots, (a_n, b_n))$ for some $n \geq 1$, where each $a_j \in S$ and $b_j$ is either 1 or 2. The ancestor process $\widehat{\xi}_s^{(x,t)}, 0 \leq s \leq t$, is a Markov process defined as follows. First, put $\widehat{\xi}_0^{(x,t)} = ((x,1))$. Suppose now that $s < t$ and $\widehat{\xi}_u^{(x,t)}$ has been defined for $u \in [0,s]$. If $\widehat{\xi}_s^{(x,t)} = \varnothing$, put $\widehat{\xi}_v^{(x,t)} = \varnothing$ for all $s < v \le t$. If $\widehat{\xi}_s^{(x,t)} = ((a_1,b_1),\dots,(a_n,b_n))$ for some $n \ge 1$, let $u$ be the the smallest time larger than $s$ at which an event occurs at (forward) time $t-u$ affecting any of the $a_j$. If there is no such $u \le t$, put $\widehat{\xi}_v^{(x,t)} = \widehat{\xi}_s^{(x,t)}$ for all $s < v \le t$, then we are done. Now suppose $u < t$:

1. If the event affecting $a_j$ at time $t-u$ is an arrow pointing from some site $a$ to $a_j$, insert $(a,b)$ into $((a_1,b_1),\dots,(a_n,b_n))$ after each $(a_i,b_i)$ such that $a_i = a_j$:

- If the arrow from $a$ to $a_j$ is labeled 2, set $b=2$ for each of these insertions.
- If the arrow is unlabeled, set $b=1$, except for insertions after any $(a_i,b_i)$ with $b_i = 2$, in which case set $b=2$.

2. If the event affecting $a_j$ is a $\delta$, then delete each $(a_i,b_i)$ from $((a_1,b_1),\dots,(a_n,b_n))$ with $a_i = a_j$.

Let $\widehat{\xi}_u^{(x,t)}$ be the resulting sequence, setting $\widehat{\xi}_u^{(x,t)} = \varnothing$ if all the $a_i$ were deleted. Iteration of this procedure defines $\widehat{\xi}_v^{(x,t)}$ for all $v \in [0,t]$. For an ancestor configuration $\widehat{\xi} = ((a_1,b_1),\dots,(a_n,b_n))$, let $\mathrm{supp}(\widehat{\xi}) = \{a_j : 1 \le j \le n\}$, so that $\mathrm{supp}(\widehat{\xi}_s^{(x,t)}) = D_s^{(x,t)}$. We say that the $j$th ancestor $a_j$ is 1-*blocked* if $b_j = 2$. If $\lambda_1 = \lambda_2$, then there are no 1-blocked ancestors and the $b_j$ can be dispensed with.

The duality equation relating $\widehat{\xi}_s^{(x,t)}$ and $\xi_{t-s}$ is

$$\xi_t(x) = \Psi(x, \widehat{\xi}_s^{(x,t)}, \xi_{t-s}), \qquad s \in [0,t], \tag{2.2}$$

where $\Psi(x,\widehat{\xi},\xi)$ is the function of sites $x \in S$, ancestor configurations $\widehat{\xi}$, and configurations $\xi \in \{0,1,2\}^S$ defined as follows. If $\widehat{\xi} = \varnothing$, put $\Psi(x,\widehat{\xi},\xi) = 0$. Otherwise, $\widehat{\xi} = ((a_1,b_1),\dots,(a_n,b_n))$ for some $n \ge 1$, and we start checking the ancestors one at a time. If $\xi(a_1) = 2$, set $\Psi(x,\widehat{\xi},\xi) = 2$, indicating a 2 propagates up. If $\xi(a_1) = 1$ and $b_1 = 1$, set $\Psi(x,\widehat{\xi},\xi) = 1$, indicating a 1 propagates up. Now suppose either $\xi(a_1) = 1$ and $b_1 = 2$, or $\xi(a_1) = 0$. If $n = 1$, set $\Psi(x,\widehat{\xi},\xi) = 0$. Otherwise $n \ge 2$, and we consider $(a_2,b_2)$ and proceed as with $(a_1,b_1)$, either setting $\Psi(x,\widehat{\xi},\xi)$ equal to 1 or 2, or exhausting the set of ancestors completely, in which case we set $\Psi(x,\widehat{\xi},\xi) = 0$. The duality equation (2.2) holds because it holds at time $s = 0$, and each transition preserves its validity. For an example, see Figure 1, in which the solid circles indicate deaths, $\xi_0(a) = 2$, $\xi_0(b) = 1$, $\xi_0(c) = 0$, $\xi_0(d) = 1$, $\xi_0(e) = 2$, $\widehat{\xi}_s^{(c,t)} = ((b,1),(a,2),(c,1),(d,1),(e,2),(d,2))$, $\xi_{t-s}(a) = \xi_{t-s}(c) = \xi_{t-s}(e) = 1, \xi_{t-s}(b) = \xi_{t-s}(d) = 0$ and $\Psi(c, \widehat{\xi}_s^{(c,t)}, \xi_{t-s}) = 1$.

The reverse time ancestor processes $\widehat{\xi}_s^{(x,t)}$ are defined only for bounded time intervals. However, as in [8], we can switch to forward time, and define ancestor processes $\bar{\xi}_s^{(x,t)}, s \geq t$, in the analogous way, so that the laws of $\bar{\xi}_s^{(x,0)}, s \in [0,t]$, and $\widehat{\xi}_s^{(x,t)}, s \in [0,t]$, are the same. Thus, in Figure 1, $\bar{\xi}_t^{(d,t-s)} = ((d,1),(e,2))$. We will write $\bar{\xi}_s^x$ for $\bar{\xi}_s^{(x,0)}$.

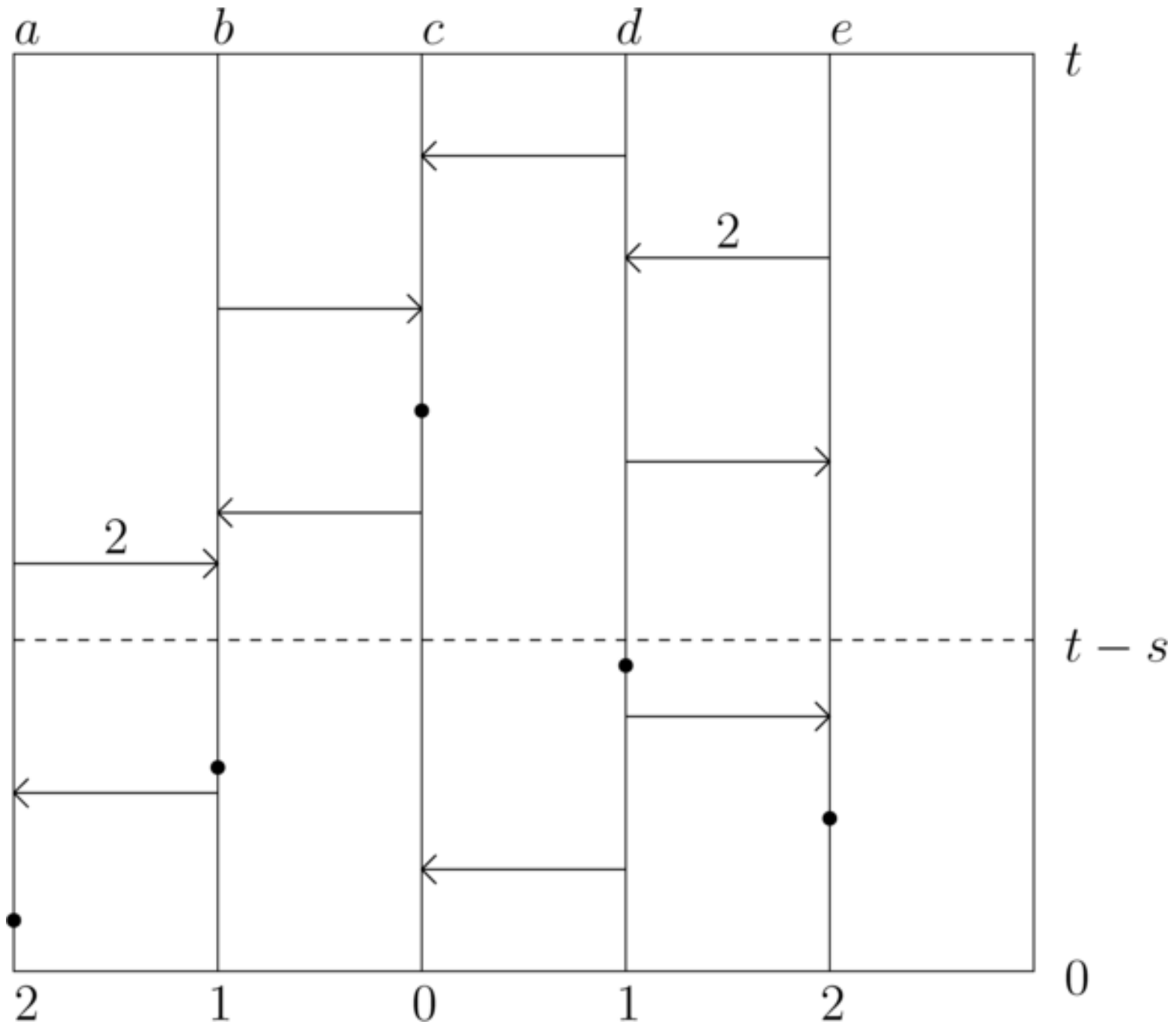

FIG. 1.

Our construction also contains various couplings with the single-type contact process. Consider $x \in S$ and $t > 0$. For $s \geq t$, let $\zeta_s^{1,(x,t)}$ be the set of sites $y$ such that there is a 1-path up from $(x,t)$ to $(y,s)$, and let $\zeta_s^{2,(x,t)}$ be the set of sites $y$ such that there is either a 1-path or 2-path up from $(x,t)$ to $(y,s)$. Write $\zeta_s^{i,x}$ for $\zeta_s^{i,(x,0)}$ and $\zeta_s^{i,A}$ for $\bigcup_{x \in A} \zeta_s^{i,x}$, $i = 1,2$. Each process $\zeta_s^{i,x}$, $s \geq 0$, is a single-type contact processes with birth rate $\lambda_i$. Also, $D_s^{1,(x,t)}, 0 \leq s \leq t$, has the same law as $\zeta_s^{1,x}, 0 \leq s \leq t$, and $D_s^{(x,t)}, 0 \leq s \leq t$, has the same law as $\zeta_s^{2,x}, 0 \leq s \leq t$. Finally, we make the observation that $\mathrm{supp}(\bar{\xi}_s^{(x,t)}), s \geq t$, and $\zeta_s^{2,(x,t)}, s \geq t$, have the same law.

A word on notation. Throughout we will use $\xi, \eta$ for elements of $\{0,1,2\}^S$, $\zeta$ for subsets of $S$, and $\bar{\xi}$ and $\widehat{\xi}$ for ancestor configurations.

## 3. Proof of Theorem 1.

Recall that (2.1) is in force.

Proof of (1.5). It suffices to assume that $\xi \equiv 1$ on $S(y)$ for some neighbor $y$ of $\mathcal{O}$, and show that there is some $x \in S(y)$ such that

$$\liminf_{t\to\infty} P(\xi_t(x) = 1) > 0. \tag{3.1}$$

We start with a simple consequence of the duality equation (2.2): for any $x \in S$, then

$$\{\varnothing \neq D_t^{1,(x,t)} \text{ and } D_t^{(x,t)} \subset S(y)\} \subset \{\xi_t(x) = 1\}. \tag{3.2}$$

This follows because on the event on the left-hand side $D_t^{(x,t)} = \mathrm{supp}(\widehat{\xi}_t^{(x,t)}) = \{a_1, \ldots, a_n\}$ for some $n \geq 1$, $\xi_0(a_i) = 1$ for each $i$, and since $D_t^{1,(x,t)}$ is nonempty, at least one of the $a_i$ is not 1-blocked. Therefore, $\Psi(x, \widehat{\xi}_t^{(x,t)}, \xi_0) = 1$, so $\xi_t(x) = 1$ by (2.2).

By (3.2), we have

$$\begin{aligned}
P(\xi_t(x) = 1) &\geq P(D_t^{1,(x,t)} \neq \varnothing) - P(D_t^{(x,t)} \not\subset S(y)) \\
&= P(\zeta_t^{1,x} \neq \varnothing) - P(\zeta_t^{2,x} \not\subset S(y)) \\
&\geq P(\zeta_s^{1,x} \neq \varnothing \text{ for all } s \geq 0) - P(\zeta_s^{2,x} \cap S^c(y) \neq \varnothing \text{ for some } s \geq 0).
\end{aligned}$$

Since $\lambda_1 > \lambda_*$, the survival probability $\alpha(\lambda_1) = \alpha_{\{x\}}(\lambda_1)$ [recall (1.8)] is positive. For $x \in S(y)$, the fact that 2's spread only by nearest neighbor contact implies

$$\{\zeta_s^{2,x} \cap S^c(y) \neq \varnothing \text{ for some } s \geq 0\} = \{\mathcal{O} \in \zeta_s^{2,x} \text{ for some } s \geq 0\}.$$

We can now make use of Theorem I.4.65 of [7], which states that if $\lambda_2 \leq \lambda^*$, then

$$P(\mathcal{O} \in \zeta_s^{2,x} \text{ for some } s \geq 0) \leq (1/\sqrt{d})^{|x-\mathcal{O}|}. \tag{3.3}$$

(All that is important for us about this bound is that it tends to 0 as $|x - \mathcal{O}| \to \infty$.) By combining the above, it follows that, for all $x \in S(y)$ and $t \geq 0$,

$$P(\xi_t(x) = 1) \geq \alpha(\lambda_1) - (1/\sqrt{d})^{|x-\mathcal{O}|}, \tag{3.4}$$

which is certainly positive for $x \in S(y)$ with $|x - \mathcal{O}|$ sufficiently large. $\square$

Before continuing with the proof of Theorem 1, we prove a lemma which shows that limits of certain finite dimensional distributions for $\xi_t$ exist when the dual lands in a region where the initial state is constant. For $x \in \mathbb{T}_d$ and $t > 0$, let $\tau_K^{(x,t)} = \inf\{s \in [0,t] : D_s^{(x,t)} \cap B_K \neq \varnothing\}$, with $\inf(\varnothing) = \infty$ [so $\tau_K^{(x,t)} = \infty$ simply means $D_s^{(x,t)} \cap B_K = \varnothing\ \forall s \leq t$]. Switching to forward time, let $\tau_K^{i,x} = \inf\{s \geq 0 : \zeta_s^{i,x} \cap B_K \neq \varnothing\}$ and define the lifetimes $\sigma^{i,x} = \inf\{s \geq 0 : \zeta_s^{i,x} = \varnothing\}$, $i = 1, 2$.

LEMMA 1. *Assume that $y \in \partial B_K$, $\xi_0$ is constant on $S(y)$, and $A_0, A_1, A_2$ are finite disjoint subsets of $S(y)$. Then as $t \to \infty$,*

$$P(\tau_K^{(z,t)} = \infty \text{ and } \xi_t(z) = j \ \forall z \in A_j, j = 0, 1, 2) \to \phi(A_0, A_1, A_2) \tag{3.5}$$

*for some $\phi(A_0, A_1, A_2)$.*

PROOF. Let $A = A_0 \cup A_1 \cup A_2 \subset S(y)$. (i) Suppose $\xi_0 \equiv 0$ on $S(y)$. For $z \in A \subset S(y)$, if $\tau_K^{(z,t)} = \infty$, then $\operatorname{supp}(\widehat{\xi}_t^{(z,t)}) \subset S(y)$. Consequently, no ancestor can land on a site occupied by a 1 or 2, implying $\xi_t(z)$ cannot be 1 or 2. This means the left-hand side of (3.5) is zero if either $A_1$ or $A_2$ is nonempty. If $A_1 = A_2 = \varnothing$, the left-hand side of (3.5) equals

$$P(\tau_K^{(z,t)} = \infty \ \forall z \in A) = P(\tau_K^{2,z} > t \ \forall z \in A) \to P(\tau_K^{2,z} = \infty \ \forall z \in A)$$

as $t \to \infty$.

(ii) Suppose $\xi_0 \equiv 2$ on $S(y)$. Again, $\tau_K^{(z,t)} = \infty$ implies $\operatorname{supp}(\widehat{\xi}_t^{(z,t)}) \subset S(y)$, so now the left-hand side of (3.5) is zero unless $A_1 = \varnothing$. In this case, $\xi_t(z) = 2$ if and only if $D_t^{(z,t)} \neq \varnothing$. Therefore, the left-hand side of (3.5), for $A_1 = \varnothing$, equals

$$\begin{aligned}
&P(\tau_K^{(z,t)} = \infty \ \forall z \in A, D_t^{(z,t)} = \varnothing \ \forall z \in A_0, D_t^{(z,t)} \neq \varnothing \ \forall z \in A_2) \\
&\quad = P(\tau_K^{2,z} > t \ \forall z \in A, \sigma^{2,z} \leq t \ \forall z \in A_0, \sigma^{2,z} > t \ \forall z \in A_2) \\
&\quad \to P(\tau_K^{2,z} = \infty \ \forall z \in A, \sigma^{2,z} < \infty \ \forall z \in A_0, \sigma^{2,z} = \infty \ \forall z \in A_2)
\end{aligned}$$

as $t \to \infty$.

(iii) Suppose $\xi_0 \equiv 1$ on $S(y)$. Now the left-hand side of (3.5) is zero unless $A_2 = \varnothing$, and in this case $\xi_t(z) = 1$ if and only if $D_t^{1,(z,t)} \neq \varnothing$. Thus, if $A_2 = \varnothing$, the left-hand side of (3.5) equals

$$\begin{aligned}
&P(\tau_K^{(z,t)} = \infty \ \forall z \in A, D_t^{1,(z,t)} = \varnothing \ \forall z \in A_0, D_t^{1,(z,t)} \neq \varnothing \ \forall z \in A_1) \\
&\quad = P(\tau_K^{2,z} > t \ \forall z \in A, \sigma^{1,z} \leq t \ \forall z \in A_0, \sigma^{1,z} > t \ \forall z \in A_1) \\
&\quad \to P(\tau_K^{2,z} = \infty \ \forall z \in A, \sigma^{1,z} < \infty \ \forall z \in A_0, \sigma^{1,z} = \infty \ \forall z \in A_1)
\end{aligned}$$

as $t \to \infty$. This proves the lemma. □

PROOF OF THEOREM 1(ii). Assume $\xi_0 \in \Xi_0$, and $K$ is such that $\xi_0$ is constant on each $S(y), y \in \partial B_K$. We will prove that, for any $x \in \mathbb{T}_d$, $a \in \{0, 1, 2\}$,

$$\lim_{t \to \infty} P(\xi_t(x) = a) \qquad \text{exists.} \tag{3.6}$$

Afterward we will show how to modify the proof to handle convergence of all finite dimensional distributions.

Let $L^{(x,t)}$ be the last time up to time $t$ the dual $D_s^{(x,t)}$ starting from $(x,t)$ contains some point of $B_K$, $L^{(x,t)} = \sup\{s \le t : B_K \cap D_s^{(x,t)} \neq \varnothing\}$, with $\sup(\varnothing) = 0$. We will prove (3.6) using a decomposition based on the value of $L^{(x,t)}$. The case $L^{(x,t)} = 0$ is easily handled. If $x \in B_K$, then necessarily $L^{(x,t)} > 0$, while for $x \notin B_K$, $\{L^{(x,t)} = 0\} = \{\tau_K^{(x,t)} = \infty\}$, and so Lemma 1 implies

$$\lim_{t\to\infty} P(\xi_t(x) = a, L^{(x,t)} = 0) \qquad \text{exists.} \tag{3.7}$$

Now suppose $L^{(x,t)} = s$ for some $s \in (0,t)$, which occurs when:

1. $\widehat{\xi}_{s-}^{(x,t)} = \widehat{\xi}'$ for some $\widehat{\xi}'$ containing a single $w \in B_K$ as one or more of the ancestors of $(x,t)$,
2. there is a $\delta$ at this site $w$ at (forward) time $t-s$,
3. $\widehat{\xi}_s^{(x,t)} = \widehat{\xi}$, where $\widehat{\xi}$ is obtained from $\widehat{\xi}'$ by removing each $(a_i, b_i)$ with $a_i = w$,
4. $D_u^{(z,t-s)} \cap B_K = \varnothing$ for all $z \in \operatorname{supp}(\widehat{\xi})$ and $u \in [0, t-s]$.

Using the convention above for $\widehat{\xi}'$ and $\widehat{\xi}$, the duality equation (2.2) implies that if $0 < t_0 < t$, we have

$$\begin{aligned} &P(\xi_t(x) = a, 0 < L^{(x,t)} < t_0) \\ &\quad = \int_0^{t_0} P(L^{(x,t)} \in ds, \xi_t(x) = a) \\ &\quad = \sum_{\widehat{\xi}'}^{(K)} \int_0^{t_0} \int P(L^{(x,t)} \in ds, \widehat{\xi}_{s-}^{(x,t)} = \widehat{\xi}', \xi_{t-s} \in d\xi) 1\{\Psi(x, \widehat{\xi}, \xi) = a\}, \end{aligned}$$

where the sum $\sum^{(K)}$ is over $\widehat{\xi}'$ such that $\operatorname{supp}(\widehat{\xi}') \cap B_K$ is a single site. By independence of disjoint space–time regions of the Poisson processes used in the construction in Section 2, and the fact that $\delta$'s occur at rate one, the above equals

$$\begin{aligned} &\sum_{\widehat{\xi}'}^{(K)} \int_0^{t_0} \int P(\widehat{\xi}_s^{(x,t)} = \widehat{\xi}') \\ &\qquad \times P(\tau_K^{(z,t-s)} = \infty \ \forall z \in \operatorname{supp}(\widehat{\xi}), \xi_{t-s} \in d\xi) 1\{\Psi(x, \widehat{\xi}, \xi) = a\}\, ds. \end{aligned}$$

For fixed $x$ and $\widehat{\xi}$, $\Psi(x, \widehat{\xi}, \xi)$ depends on $\xi$ only through the values $\xi(z)$, $z \in \operatorname{supp}(\widehat{\xi})$. Hence, we can define $\Pi(x, a, \widehat{\xi})$ to be the set of partitions $\mathcal{A} = (A_0, A_1, A_2)$ of $\operatorname{supp}(\widehat{\xi})$ such that, for all $\xi$, $\Psi(x, \widehat{\xi}, \xi) = a$ if and only if $\xi \equiv j$ on $A_j, j = 0, 1, 2,$ for some $\mathcal{A} \in \Pi(x, a, \widehat{\xi})$. Therefore, summing over the possible values of $\xi_{t-s}$ on $\operatorname{supp}(\widehat{\xi})$, and replacing $P(\widehat{\xi}_s^{(x,t)} = \widehat{\xi}')$ with $P(\bar{\xi}_s^x = \widehat{\xi}')$, we have

$$\begin{aligned} &P(\xi_t(x) = a, 0 < L^{(x,t)} \le t_0) \\ &\quad = \sum_{\widehat{\xi}'}^{(K)} \sum_{\mathcal{A}\in\Pi(x,a,\widehat{\xi})} \int_0^{t_0} P(\bar{\xi}_s^x = \widehat{\xi}') q_t(s,\mathcal{A})\, ds, \end{aligned} \tag{3.8}$$

where

$$q_t(s,\mathcal{A}) = P(\tau_K^{(z,t-s)} = \infty \text{ and } \xi_{t-s}(z) = j\ \forall z \in A_j, j = 0,1,2). \tag{3.9}$$

It is time to use the fact that $\xi_0 \in \Xi_0$. Let $S^i$ be the union of all $S(y), y \in \partial B_K$ such that $\xi_0 \equiv i$ on $S(y)$, and for the sets $A_j$ above let $A_j^i = A_j \cap S^i$. Since the $S^i$ are disjoint, and the duals $D_u^{(z,t-s)}$ do not leave their respective sectors $S^i$ on $\{\tau_K^{(z,t-s)} = \infty\}$, independence of disjoint space–time regions implies

$$\begin{aligned} q_t(s,\mathcal{A}) &= P\left(\bigcap_{i=0}^{2} \{\tau_K^{(z,t-s)} = \infty \text{ and } \xi_{t-s}(z) = j\ \forall z \in A_j^i, j = 0,1,2\}\right) \\ &= \prod_{i=0}^{2} P(\tau_K^{(z,t-s)} = \infty \text{ and } \xi_{t-s}(z) = j\ \forall z \in A_j^i, j = 0,1,2). \end{aligned}$$

By Lemma 1, the above product converges for fixed $s$ as $t \to \infty$ to some $\phi(\mathcal{A}) = \prod_{i=0}^2 \phi(A_0^i, A_1^i, A_2^i)$. Since each $\Pi(x,a,\widehat{\xi}')$ is finite, this implies

$$\lim_{t\to\infty} \sum_{\mathcal{A}\in\Pi(x,a,\widehat{\xi})} q_t(s,\mathcal{A}) = \sum_{\mathcal{A}\in\Pi(x,a,\widehat{\xi})} \phi(\mathcal{A}) \tag{3.10}$$

for fixed $s, \widehat{\xi}'$.

Since $\sum_{\mathcal{A}\in\Pi(x,a,\widehat{\xi})} q_t(s,\mathcal{A}) \le 1$ and $\sum_{\widehat{\xi}'}^{(K)} \int_0^{t_0} P(\bar{\xi}_s^x = \widehat{\xi}')\, ds \le t_0$, the dominated convergence theorem can be applied in (3.8), so that (3.10) implies

$$\begin{aligned} &\lim_{t\to\infty} P(\xi_t(x) = a, 0 < L^{(x,t)} \le t_0) \\ &\quad = \sum_{\widehat{\xi}'}^{(K)} \sum_{\mathcal{A}\in\Pi(x,a,\widehat{\xi})} \phi(\mathcal{A}) \int_0^{t_0} P(\bar{\xi}_s^x = \widehat{\xi}')\, ds. \end{aligned} \tag{3.11}$$

The contribution to (3.6) of $L^{(x,t)} > t_0$ for large $t_0$ is negligible, since $\lambda_2 \le \lambda^*$ implies

$$\begin{aligned} &\limsup_{t_0\to\infty} \sup_{t\ge t_0} P(t_0 < L^{(x,t)} \le t) \\ &\quad \le \limsup_{t_0\to\infty} P(\zeta_s^{2,x} \cap B_K \ne \varnothing \text{ for some } s \ge t_0) = 0. \end{aligned} \tag{3.12}$$

We have thus proved

$$\lim_{t\to\infty} P(\xi_t(x)=a, 0<L^{(x,t)}\le t) = \sum_{\widehat{\xi}'}^{(K)} \sum_{\mathcal{A}\in\Pi(x,a,\widehat{\xi})} \phi(\mathcal{A}) \int_0^\infty P(\bar{\xi}_s^x = \widehat{\xi}')\,ds$$

and, in view of (3.7), this implies (3.6) must hold.

More generally, let $\Gamma=\{x_1,\dots,x_m\}\subset\mathbb{T}_d$, and $a:\Gamma\to\{0,1,2\}$. We claim that

$$\lim_{t\to\infty} P(\xi_t(x)=a(x)\ \forall x\in\Gamma) \qquad \text{exists.} \tag{3.13}$$

First, let $L^{\Gamma,t}=\max\{L^{(x,t)}, x\in\Gamma\}$, and use Lemma 1 to obtain

$$\lim_{t\to\infty} P(\xi_t(x)=a(x)\ \forall x\in\Gamma, L^{\Gamma,t}=0) \qquad \text{exists.}$$

Next, consider $P(\xi_t(x)=a(x)\ \forall x\in\Gamma, 0<L^{\Gamma,t}\le t_0)$, and decompose this event according to the value of $L^{\Gamma,t}$,

$$\begin{aligned} &P(\xi_t(x)=a(x)\ \forall x\in\Gamma, 0<L^{\Gamma,t}\le t_0)\\ &\quad= \sum_{\widehat{\xi}'_1,\dots,\widehat{\xi}'_m}^{(K)} \int_0^{t_0}\int P(L^{\Gamma,t}\in ds, \widehat{\xi}_{s-}^{(x_i,t)}=\widehat{\xi}'_i, 1\le i\le m, \xi_{t-s}\in d\xi)\\ &\qquad\times\prod_{i=1}^m 1\{\Psi(x_i,\widehat{\xi}_i,\xi)=a(x_i)\}, \end{aligned}$$

where the sum $\sum^{(K)}$ is over $\widehat{\xi}'_1,\dots,\widehat{\xi}'_m$ such that $B_K\cap(\bigcup_i \operatorname{supp}(\widehat{\xi}'_i))$ is a single site $w\in B_K$. As before, for each $i$, $\widehat{\xi}_i$ is obtained from $\widehat{\xi}'_i$ by removing all $(a,b)$ with $a=w$, if any. By independence of disjoint space–time regions, the above equals

$$\begin{aligned} &\sum_{\widehat{\xi}'_1,\dots,\widehat{\xi}'_m}^{(K)} \int_0^{t_0}\int P(\widehat{\xi}_s^{(x_i,t)}=\widehat{\xi}'_i, 1\le i\le m)\\ &\qquad\qquad\times P\left(\tau_K^{(z,t-s)}=\infty\ \forall z\in\operatorname{supp}\left(\bigcup_i\widehat{\xi}_i\right), \xi_{t-s}\in d\xi\right)\\ &\qquad\quad\times\prod_{i=1}^m 1\{\Psi(x_i,\widehat{\xi}_i,\xi)=a(x_i)\}\,ds. \end{aligned}$$

Let $\Pi(\Gamma,a,\widehat{\xi}_1,\dots,\widehat{\xi}_m)$ be the set of partitions $\mathcal{A}=(A_0,A_1,A_2)$ of $\bigcup_i\operatorname{supp}(\widehat{\xi}_i)$ such that $\Psi(x_i,\widehat{\xi}_i,\xi)=a(x_i)\ \forall 1\le i\le m$ if and only if $\xi\equiv j$ on $A_j$, $j=0,1,2$, for some $\mathcal{A} \in \Pi(\Gamma, a, \widehat{\xi}_1, \ldots, \widehat{\xi}_m)$. Reasoning as before gives

$$\begin{aligned} &P(\xi_t(x) = a(x)\ \forall x \in \Gamma, 0 < L^{\Gamma,t} < t_0) \\ &\quad = \sum_{\widehat{\xi}'_1,\ldots,\widehat{\xi}'_m}^{(K)} \sum_{\mathcal{A}\in\Pi(\Gamma,a,\widehat{\xi}'_1,\ldots,\widehat{\xi}'_m)} \int_0^{t_0} P(\bar{\xi}_s^{x_i} = \widehat{\xi}'_i, 1 \le i \le m) \\ &\qquad \times P(\tau_K^{(z,t-s)} = \infty \text{ and } \xi_{t-s}(z) = j\ \forall z \in A_j, j = 0,1,2)\, ds. \end{aligned}$$

By the argument leading to (3.11), for some $\phi(\mathcal{A})$,

$$\begin{aligned} &\lim_{t\to\infty} P(\xi_t(x) = a(x)\ \forall x \in \Gamma, 0 < L^{\Gamma,t} \le t_0) \\ &\quad = \sum_{\widehat{\xi}'_1,\ldots,\widehat{\xi}'_m}^{(K)} \sum_{\mathcal{A}\in\Pi(\Gamma,a,\widehat{\xi}'_1,\ldots,\widehat{\xi}'_m)} \phi(\mathcal{A}) \int_0^{t_0} P(\widehat{\xi}_s^{x_i} = \widehat{\xi}'_i, 1 \le i \le m)\, ds. \end{aligned}$$

In view of (3.12), this completes the proof of (3.13). $\square$

Proof of Theorem 1(iii). We may suppose $\lambda_1 \le \lambda_2$, $i = 1$ and $j = 2$. Our first task is to prove $\nu_\eta(\xi : |{}^1\xi| = \infty) = 1$. To do this, it is enough by (ii) to show that for $\varepsilon > 0$ and $M \ge 1$ there is a finite $A \subset S(y)$ such that

$$\lim_{t\to\infty} P\left(\sum_{x\in A} 1\{\xi_t(x) = 1\} \ge M\right) \ge 1 - \varepsilon. \tag{3.14}$$

To do this, let $\ell > 0$ be large enough so that $(1/\sqrt{d})^\ell \le \alpha(\lambda_1)/2$. Letting $\varepsilon_0 = \alpha(\lambda_1)/2 > 0$, it follows from (3.3) that, for any $y_0 \in S(y)$ and $x_0 \in S(y_0)$ such that $|x_0 - y_0| \ge \ell$,

$$P(\xi_t(x_0) = 1 \text{ and } y_0 \notin D_s^{(x_0,t)}\ \forall s \le t) \ge \varepsilon_0 \qquad \forall t \ge 0.$$

We construct the set $A$ in (3.14) as follows. Let $N$ be a positive integer large enough so that a binomial random variable $X$ with parameters $N$ and $p \ge \varepsilon_0$ will satisfy $P(X \ge M) \ge 1 - \varepsilon$. Let $y_1, \ldots, y_N$ be vertices in $S(y)$ such that $S(y_j) \cap S(y_k) = \varnothing$ for $j \ne k$, and let $x_j \in S(y_j)$ satisfy $|x_j - y_j| \ge \ell$ for $j = 1, \ldots, N$. Fix $t > 0$ and define $\varepsilon_j = 1\{\xi_t(x_j) = 1, y_j \notin D_s^{(x_j,t)}\ \forall s \le t\}$, $j = 1, \ldots, N$. By independence of disjoint space–time regions, the $\varepsilon_j$ are independent with $P(\varepsilon_j = 1) \ge \varepsilon_0$, so $X = \sum_{j=1}^N \varepsilon_j$ is binomial with parameters $N$ and $p \ge \varepsilon_0$ and thus, (3.14) holds for $A = \{x_i, i = 1, \ldots, N\}$.

The fact that $\nu(\xi : \xi(x) = 2) \to 0$ as $x \to \infty$, $x \in S(y)$ is a simple consequence of duality and the bound (3.3). Since $x \in S(y)$ and $\xi_0 \equiv 1$ on $S(y)$,

$$P(\xi_t(x) = 2) \le P(\widehat{\xi}_t^{(x,t)} \cap S^c(y) \ne \varnothing) \le (1/\sqrt{d})^{|x-y|}.$$

Now let $t \to \infty$ and $x \to \infty, x \in S(y)$. $\square$

PROOF OF THEOREM 1(iv). Assume $\eta, \eta' \in \Xi_0$. Consider both $\xi_t$ with $\xi_0 = \eta$ and $\xi'_t$ with $\xi'_0 = \eta'$ defined via the Poisson processes in Section 2. Let $A = \{y \in \mathbb{T}_d : \eta(y) \neq \eta'(y)\}$, and suppose $A$ is finite. For any $x$, since each $\lambda_i \leq \lambda^*$, $P(\widehat{\xi}_t^{(x,t)} \cap A \neq \varnothing) \to 0$ as $t \to \infty$. By the duality equation (2.2), since $\xi_0 = \xi'_0$ on $A^c$, $P(\xi_t(x) \neq \xi'_t(x)) \to 0$ as $t \to \infty$, which implies $\nu_\eta = \nu_{\eta'}$. The remaining conclusion of (iv) is a simple consequence of (iii). □

**4. Proof of Theorem 2.** Before beginning the proof of Theorem 2, we state and prove a fact about the upper invariant measure $\bar{\nu}$ for the single-type contact process. The result we need is an immediate consequence of the inequality I.2.30(b) of [7] for $S = \mathbb{Z}^d$. This inequality also holds for $S = \mathbb{T}_d$, but we could not find a reference. Since the following result is all we need, and its proof may apply to other choices of $S$, we give the proof here.

PROPOSITION 1. *Assume $S$ is $\mathbb{Z}^d$ or $\mathbb{T}_d$ and $\bar{\nu}$ is the upper invariant measure for the contact process with birth rate $\lambda > \lambda_*$. Let $\delta_L = \sup_{|A|>L} \bar{\nu}(\zeta : \zeta \cap A = \varnothing)$. Then*

$$\delta_L \to 0 \qquad \text{as } L \to \infty. \tag{4.1}$$

PROOF. Recall the survival probabilities $\alpha_A$ from (1.8). A consequence of duality for the contact process (see (I.1.8) of [7]) is that

$$\bar{\nu}(\zeta : \zeta \cap A \neq \varnothing) = \alpha_A. \tag{4.2}$$

Let $\zeta_t$ denote the single-type contact process with birth rate $\lambda$. Define the *lifetime* $\sigma(x)$ and the *radius* $r(x)$ of the process started at $x$ by $\sigma(x) = \inf\{t \geq 0 : \zeta_t^x = \varnothing\}$ and $r(x) = \sup\{|y - x| : y \in \zeta_t^x \text{ for some } t \geq 0\}$. Let $\varepsilon_M(x) = P(\sigma(x) < \infty, r(x) \geq M)$. By translation invariance, $\varepsilon_M = \varepsilon_M(x)$ does not depend on $x$, and $\varepsilon_M \to 0$ as $M \to \infty$.

Now fix $M, N > 0$ and let $A$ be any set consisting of $N$ points such that $|x - y| \geq M$ for all $x, y \in A, x \neq y$. Then, by the construction in Section 2,

$$\begin{aligned} 1 - \alpha_A &= P(\sigma(x) < \infty \ \forall x \in A) \\ &\leq P(\sigma(x) < \infty, r(x) < M \ \forall x \in A) \\ &\quad + P(\exists x \in A : \sigma(x) < \infty \text{ and } r(x) \geq M). \end{aligned}$$

By independence of disjoint space–time regions and translation invariance, the first probability on the right-hand side equals $\prod_{x \in A} P(\sigma(x) < \infty, r(x) < M) \leq (1 - \alpha)^N$ and the second probability is bounded above by $\sum_{x \in A} \varepsilon_M(x) \leq N\varepsilon_M$. Thus,

$$1 - \alpha_A \leq (1 - \alpha)^N + N\varepsilon_M.$$

Since $\alpha > 0$, given $\varepsilon > 0$, we may choose first $N$ and then $M$ so that $(1 - \alpha)^N < \varepsilon/2$ and $N\varepsilon_M < \varepsilon/2$. Now let $L_0$ be large enough so that any $A \subset S$ with at least $L_0$ points must contain at least $N$ points separated from one another by distance at least $M$. It follows from the monotonicity of $\alpha_A$ that if $|A| > L_0$, then $\bar{\nu}(\zeta : \zeta \cap A = \varnothing) = 1 - \alpha_A < \varepsilon$. $\square$

Proof of Theorem 2. Here is the idea of the proof. Suppose $t, u > 0$ are large, $\xi_{t+u}(x) \neq 0$ and ${}^2\xi_{t+u} \neq \varnothing$. Looking forward from time 0, ${}^2\xi_u$ cannot be empty and, with high probability, will have many points. Looking backward from time $t+u$, the dual ancestor process starting at $(x, t+u)$ must survive $t$ time units. For $T > 0$ large but small compared to $t$, we search for a space–time point $(y, t+u-s)$ such that $y \in D_s^{(x,t+u)} = \operatorname{supp}(\widehat{\xi}_s^{(x,t+u)})$ is the primary ancestor at time $s$, $y$ is 1-blocked, and $D_T^{(y,t+u-s)}$ is nonempty. Trying at most a geometric number of times, with high probability, we will find such a point $(y, s)$ with $s$ not too large. Furthermore, $D_T^{(y,t+u-s)} \neq \varnothing$ will imply that $D_{t-s}^{(y,t+u-s)} \neq \varnothing$ with high probability, and also that $D_{t-s}^{(y,t+u-s)}$ will intersect ${}^2\xi_u$. (It is this last point which fails unless $\lambda_2 > \lambda^*$.) This will prevent $\xi_{t+u}(x) = 1$ since the sites in $D_{t-s}^{(y,t+u-s)}$ are descendants of a 1-blocked primary ancestor.

We prepare for the proof of (1.7) by assembling a few preliminary facts. Recall $\bar{\nu} = \bar{\nu}_{\lambda_2}$ from (1.9). Since $\lambda_2 > \lambda^*$, the complete convergence theorem (1.10) implies that, for any site $x$ and finite $A \subset S$, with $\alpha = \alpha(\lambda_2)$,

$$\lim_{t\to\infty} P(\zeta_t^{2,x} \neq \varnothing, \zeta_t^{2,x} \cap A = \varnothing) = \alpha\bar{\nu}(\zeta : \zeta \cap A = \varnothing). \tag{4.3}$$

Let $\rho(T) = P(\zeta_T^{2,x} \neq \varnothing, \zeta_t^{2,x} = \varnothing \text{ for some } t > T)$, and observe that

$$\rho(T) \to 0 \qquad \text{as } T \to \infty. \tag{4.4}$$

Since individuals die at constant rate one, standard arguments show that the 2's must either die out or that their number must tend to infinity, that is,

$$P(1 \le |{}^2\xi_u| \le L) \to 0 \qquad \text{as } u \to \infty \tag{4.5}$$

for fixed $L > 0$.

Our argument uses a certain subset $A_t^x$ of the highest priority ancestors of the forward time ancestor process $\bar{\xi}_t^x$ which we now define. If $\bar{\xi}_t^x = \varnothing$, put $A_t^x = \varnothing$. Now suppose that $\bar{\xi}_t^x = ((a_1(t), b_1(t)), \ldots, (a_n(t), b_n(t))$ for some $n \ge 1$. Put $A_t^x = \varnothing$ if $a_1(t)$ is not 1-blocked [i.e., $b_1(t) = 1$], and otherwise let

$$A_t^x = \{a_1(t), \ldots, a_m(t)\}, \tag{4.6}$$

where $m$ is the largest index such that $a_1(t), \ldots, a_m(t)$ are all 1-blocked. We will see below that, for large $t$, $A_t^x$ is large whenever $\bar{\xi}_t^x \neq \varnothing$.

Our goal is to prove

$$\lim_{u\to\infty} \limsup_{t\to\infty} P(\xi_{t+u}(x) = 1, |{}^2\xi_u| \ge 1) = 0. \tag{4.7}$$

This implies (1.7) since $\{\xi_{t+u}(x)=1, |^2\xi_{t+u}|\geq 1\} \subset \{\xi_{t+u}(x)=1, |^2\xi_u|\geq 1\}$. In fact, on account of (4.5), we may focus our attention on $\{\xi_{t+u}(x) = 1, |^2\xi_u| > L\}$ for large $L$.

We begin with an application of the duality equation (2.2),

$$P(\xi_{t+u}(x)=1, |^2\xi_u|>L) = P(\Psi(x, \widehat{\xi}_t^{(x,t+u)}, \xi_u)=1, |^2\xi_u|>L). \tag{4.8}$$

By independence of disjoint space–time regions, and switching to the forward time ancestor process $\bar{\xi}_t^x$ in the second equality,

$$\begin{aligned} &P(\Psi(x, \widehat{\xi}_t^{(x,t+u)}, \xi_u)=1, |^2\xi_u|>L) \\ &\qquad = \int_{|^2\eta|>L} P(\xi_u\in d\eta) P(\Psi(x, \widehat{\xi}_t^{(x,t+u)}, \eta)=1) \\ &\qquad = \int_{|^2\eta|>L} P(\xi_u\in d\eta) P(\Psi(x, \bar{\xi}_t^{x}, \eta)=1). \end{aligned} \tag{4.9}$$

It follows from the the definitions of $\Psi$ and $A_t^x$ that

$$\Psi(x, \bar{\xi}_t^x, \eta) = 2 \qquad \text{on } \{\bar{\xi}_t^x \neq \varnothing, A_t^x \cap {}^2\eta \neq \varnothing\}, \tag{4.10}$$

which implies

$$P(\Psi(x, \bar{\xi}_t^x, \eta)=1) \leq P(\bar{\xi}_t^x \neq \varnothing, A_t^x\cap {}^2\eta = \varnothing). \tag{4.11}$$

Now combine this with (4.8) and (4.9) to obtain

$$\begin{aligned} &P(\xi_{t+u}(x)=1, |^2\xi_u|>L) \\ &\qquad \leq \int_{|^2\eta|>L} P(\xi_u\in d\eta) P(\bar{\xi}_t^x \neq \varnothing, A_t^x\cap {}^2\eta = \varnothing). \end{aligned} \tag{4.12}$$

To show the right-hand side above is small, we will argue that, for large $t$ and $u$, $A_t^x \cap {}^2\eta = \varnothing$ is unlikely when $\bar{\xi}_t^x \neq \varnothing$ by means of the following construction. For each $y\in S$ pick some nearest neighbor $\tilde{y}$, and fix this assignment. Call a space–time point $(y,s)$ *good* if the following events happen:

1. a 2-only arrow pointing from $\tilde{y}$ to $y$ occurs during $(s, s+1)$,
2. a $\delta$ occurs at $y$ at some time during $(s+1, s+2)$,
3. no other events affecting $y$ or $\tilde{y}$ occur during $[s, s+2]$.

Then $\varepsilon_0 = P((y,s) \text{ is good}) > 0$ and does not depend on $(y,s)$. If $\bar{\xi}_s^x \neq \varnothing$, then $\bar{\xi}_s^x = ((a_1(s), b_1(s)), \dots, (a_n(s), b_n(s)))$ for some $n\geq 1$ and primary ancestor site $a_1(s)$. If $(a_1(s), s)$ is good, then $a_1(s+2) = a_1(s)$, and for $v\geq s+2$ the sites of $\bar{\xi}_v^{(a_1(s+2), s+2)}$ are the highest priority sites of $\bar{\xi}_v^x$, and they are all 1-blocked.

Fix $T>0$ and let $s_k=k(T+2)$ and $t_k=s_k+2$, $k\ge 0$. Let $R$ be the smallest $k$ such that the primary ancestor $a$ in $\bar\xi^x_{s_k}$ is good at time $s_k$ and the ancestor process starting at $(\tilde a,t_k)$ lasts at least $T$ time units, that is,

$$R=\inf\{k:\bar\xi^x_{s_k}\neq\varnothing,(a_1(s_k),s_k)\text{ is good and }\bar\xi^{(a_1(t_k),t_k)}_{s_{k+1}}\neq\varnothing\}.$$

Let $\varepsilon_1=\varepsilon_0\alpha(\lambda_2)>0$. By independence of disjoint space–time regions, for any $(y,s)$ and $T>0$,

$$\begin{aligned}&P((y,s)\text{ is good and }\bar\xi^{(\tilde y,s+2)}_{s+2+T}\neq\varnothing)\\&\qquad=P((y,s)\text{ is good})P(\bar\xi^{(\tilde y,s+2)}_{s+2+T}\neq\varnothing)\\&\qquad=P((y,s)\text{ is good})P(\zeta^{2,(\tilde y,s+2)}_{s+2+T}\neq\varnothing)\ge\varepsilon_1.\end{aligned}\tag{4.13}$$

Iterating this inequality and using the Markov property gives us

$$P(\bar\xi^x_{s_k}\neq\varnothing\text{ and }R>k)\le(1-\varepsilon_1)^k,\qquad k\ge 0.\tag{4.14}$$

Consequently, if $k_0>0$ and $t>s_{k_0+1}$,

$$\begin{aligned}&P(\bar\xi^x_t\neq\varnothing,A^x_t\cap{}^2\eta=\varnothing)\\&\qquad\le(1-\varepsilon_1)^{k_0}+\sum_{k=0}^{k_0}P(R=k,\bar\xi^x_t\neq\varnothing,A^x_t\cap{}^2\eta=\varnothing).\end{aligned}\tag{4.15}$$

Now define the events

$$\begin{aligned}G_k(a)=\{&R>k-1,\bar\xi^x_{s_k}\neq\varnothing,\\&a_1(s_k)=a,\text{ and }(a,s_k)\text{ is good}\},\qquad a\in S.\end{aligned}\tag{4.16}$$

For $k\le k_0$, $\operatorname{supp}(\bar\xi^{(\tilde a,t_k)}_t)\subset A^x_t$ on $G_k(a)$, and this implies

$$\begin{aligned}&P(R=k,\bar\xi^x_t\neq\varnothing,A^x_t\cap{}^2\eta=\varnothing)\\&\qquad\le\sum_{a\in S}P(G_k(a)\cap\{\bar\xi^{(\tilde a,t_k)}_{s_{k+1}}\neq\varnothing,\operatorname{supp}(\bar\xi^{(\tilde a,t_k)}_t)\cap{}^2\eta=\varnothing\}).\end{aligned}$$

For each $a$, by independence of disjoint space–time regions,

$$\begin{aligned}&P(G_k(a)\cap\{\bar\xi^{(\tilde a,t_k)}_{s_{k+1}}\neq\varnothing,\ \operatorname{supp}(\bar\xi^{(\tilde a,t_k)}_t)\cap{}^2\eta=\varnothing\})\\&\qquad\le P(G_k(a))P(\bar\xi^{(\tilde a,t_k)}_{s_{k+1}}\neq\varnothing,\operatorname{supp}(\bar\xi^{(\tilde a,t_k)}_t)\cap{}^2\eta=\varnothing)\\&\qquad=P(G_k(a))P(\zeta^{2,\tilde a}_T\neq\varnothing,\zeta^{2,\tilde a}_{t-t_k}\cap{}^2\eta=\varnothing)\\&\qquad\le P(G_k(a))(\rho(T)+P(\zeta^{2,\tilde a}_{t-t_k}\neq\varnothing,\zeta^{2,\tilde a}_{t-t_k}\cap{}^2\eta=\varnothing)),\end{aligned}\tag{4.17}$$

where we have shifted back to time 0 and used the fact that $t > s_{k_0+1}$. Therefore, by the bounds (4.15) and (4.17), and the fact that $\sum_{a\in S} P(G_k(a)) \le 1$,

$$\begin{aligned} &P(\bar\xi_t^x \neq \varnothing, A_t^x \cap {}^2\eta = \varnothing) \\ &\quad \le (1-\varepsilon_1)^{k_0} + \rho(T)(k_0+1) \\ &\qquad + \sum_{k=0}^{k_0}\sum_{a\in S} P(G_k(a))P(\zeta_{t-t_k}^{2,\tilde a} \neq \varnothing, \zeta_{t-t_k}^{2,\tilde a} \cap {}^2\eta = \varnothing). \end{aligned} \tag{4.18}$$

Recall the definition of $\delta_L$ in Proposition 1. For fixed $k$, $a$ and $\eta$ such that $|{}^2\eta| > L$, (4.3) implies that

$$\lim_{t\to\infty} P(\zeta_{t-t_k}^{2,\tilde a} \neq \varnothing, \zeta_{t-t_k}^{2,\tilde a} \cap {}^2\eta = \varnothing) = \alpha\bar\nu(\zeta : \zeta \cap {}^2\eta = \varnothing) \le \delta_L.$$

With this we can apply Fatou to (4.18) to obtain

$$\limsup_{t\to\infty} P(\bar\xi_t^x \neq \varnothing, A_t^x \cap {}^2\eta = \varnothing) \le (1-\varepsilon_1)^{k_0} + (k_0+1)(\rho(T)+\delta_L) \tag{4.19}$$

and letting $T\to\infty$ and using (4.4) then gives us

$$\limsup_{t\to\infty} P(\bar\xi_t^x \neq \varnothing, A_t^x \cap {}^2\eta = \varnothing) \le (1-\varepsilon_1)^{k_0} + (k_0+1)\delta_L. \tag{4.20}$$

With this inequality, we apply Fatou again, this time to (4.12), to obtain

$$\limsup_{t\to\infty} P(\xi_{t+u}(x) = 1, |{}^2\xi_u| > L) \le (1-\varepsilon_1)^{k_0} + (k_0+1)\delta_L. \tag{4.21}$$

We have finally established the bound

$$\begin{aligned} &\limsup_{t\to\infty} P(\xi_{t+u}(x) = 1, |{}^2\xi_u| \ge 1) \\ &\qquad \le P(1 \le |{}^2\xi_u| \le L) + (1-\varepsilon_1)^{k_0} + (k_0+1)\delta_L. \end{aligned} \tag{4.22}$$

Now let $u, L, k_0 \to \infty$ in order and use (4.5) and (4.1) to finish the proof of (4.7). □

**5. Proof of Theorem 3.** Recall the definitions (1.8) and (1.11) and the assumptions of Theorem 3. In view of Theorem 2, it is enough to prove that for all finite $A \subset S$,

$$P({}^i\xi_t \cap A \neq \varnothing) \to \alpha_\eta^i \bar\nu_{\lambda_i}(\zeta : \zeta \cap A \neq \varnothing) \qquad \text{as } t\to\infty, i = 1,2. \tag{5.1}$$

By (4.2), $\bar\nu_{\lambda_i}(\zeta : \zeta \cap A \neq \varnothing) = \alpha_A(\lambda_i)$ and, hence, (5.1) is equivalent to

$$\lim_{t\to\infty} P({}^i\xi_t \cap A \neq \varnothing) = \alpha_\eta^i \alpha_A(\lambda_i), \qquad i = 1,2. \tag{5.2}$$

First consider the case $i = 1$. By (1.7), (5.2) will follow once we establish

$$\lim_{t\to\infty} P({}^1\xi_t \cap A \neq \varnothing, {}^2\xi_t = \varnothing) = \alpha_\eta^1 \alpha_A(\lambda_1). \tag{5.3}$$

By independence of disjoint space–time regions, for $t, u > 0$,

$$
\begin{aligned}
&P({}^1\xi_{t+u} \cap A \neq \varnothing, {}^2\xi_u = \varnothing) \\
&\quad \leq P({}^1\xi_u \neq \varnothing, {}^2\xi_u = \varnothing, D_t^{1,(x,t+u)} \neq \varnothing \text{ for some } x \in A) \\
&\quad = P({}^1\xi_u \neq \varnothing, {}^2\xi_u = \varnothing) P(D_t^{1,(x,t+u)} \neq \varnothing \text{ for some } x \in A) \\
&\quad = P({}^1\xi_u \neq \varnothing, {}^2\xi_u = \varnothing) P(\zeta_t^{1,A} \neq \varnothing) \\
&\quad \to \alpha_\eta^1 \alpha_A(\lambda_1) \qquad \text{as } t, u \to \infty.
\end{aligned}
$$

Since

$$
\begin{aligned}
&P({}^2\xi_u \neq \varnothing, {}^2\xi_{t+u} = \varnothing) \\
&\quad \leq P({}^2\xi_u \neq \varnothing, {}^2\xi_s = \varnothing \text{ for some } s \geq u) \to 0 \qquad \text{as } u \to \infty,
\end{aligned}
\tag{5.4}
$$

this proves

$$
\limsup_{t\to\infty} P({}^1\xi_t \cap A \neq \varnothing, {}^2\xi_t = \varnothing) \leq \alpha_\eta^1 \alpha_A(\lambda_1). \tag{5.5}
$$

For the required lower bound, duality implies

$$
\begin{aligned}
&P({}^1\xi_{t+u} \cap A \neq \varnothing, {}^2\xi_{t+u} = \varnothing) \\
&\quad \geq P({}^2\xi_u = \varnothing, \Psi(x, \widehat{\xi}_t^{(x,t+u)}, \xi_u) = 1 \text{ for some } x \in A).
\end{aligned}
$$

On the event $\{{}^2\xi_u = \varnothing\}$, $\Psi(x, \widehat{\xi}_t^{(x,t+u)}, \xi_u) = 1$ if and only if $D_t^{1,(x,t+u)} \cap {}^1\xi_u \neq \varnothing$. Consequently, for any $L > 0$,

$$
\begin{aligned}
&P({}^1\xi_{t+u} \cap A \neq \varnothing, {}^2\xi_{t+u} = \varnothing) \\
&\quad \geq \int_{|{}^1\eta|>L, {}^2\eta=\varnothing} P(\xi_u \in d\eta) P(D_t^{1,(x,t+u)} \cap {}^1\eta \neq \varnothing \text{ for some } x \in A) \\
&\quad = \int_{|{}^1\eta|>L, {}^2\eta=\varnothing} P(\xi_u \in d\eta) P(\zeta_t^{1,A} \cap {}^1\eta \neq \varnothing).
\end{aligned}
$$

Now we replace $P(\zeta_t^{1,A} \cap {}^1\eta \neq \varnothing)$ with $P(\zeta_t^{1,A} \neq \varnothing) - P(\zeta_t^{1,A} \neq \varnothing, \zeta_t^{1,A} \cap {}^1\eta = \varnothing)$ so that the above implies

$$
\begin{aligned}
&P({}^1\xi_{t+u} \cap A \neq \varnothing, {}^2\xi_{t+u} = \varnothing) \\
&\quad \geq P(|{}^1\xi_u| > L, {}^2\xi_u = \varnothing) P(\zeta_t^{1,A} \neq \varnothing) \\
&\quad\quad - \int_{|{}^1\eta|>L, {}^2\eta=\varnothing} P(\xi_u \in d\eta) P(\zeta_t^{1,A} \neq \varnothing, \zeta_t^{1,A} \cap {}^1\eta = \varnothing).
\end{aligned}
\tag{5.6}
$$

For fixed $u$ and $L$,

$$
\begin{aligned}
&\lim_{t\to\infty} P(|{}^1\xi_u| > L, {}^2\xi_u = \varnothing) P(\zeta_t^{1,A} \neq \varnothing) \\
&\quad \geq P(|{}^1\xi_u| \geq 1, {}^2\xi_u = \varnothing) \alpha_A(\lambda_1) - P(1 \leq |{}^1\xi_u| \leq L).
\end{aligned}
\tag{5.7}
$$

For fixed $\eta$ with $|{}^1\eta| > L$, the complete convergence theorem (1.10) implies that

$$\lim_{t\to\infty} P(\zeta_t^{1,A} \neq \varnothing, \zeta_t^{1,A} \cap {}^1\eta = \varnothing) = \alpha_A \bar{\nu}_{\lambda_1}(\zeta : \zeta \cap \eta = \varnothing) \le \delta_L \tag{5.8}$$

(recall $\delta_L$ from Proposition 1). We can now plug (5.7) and (5.8) into (5.6) and use Fatou, keeping $u$ and $L$ fixed, to get

$$\begin{aligned} &\liminf_{t\to\infty} P({}^1\xi_{t+u} \cap A \neq \varnothing, {}^2\xi_{t+u} = \varnothing) \\ &\qquad \ge P(|{}^1\xi_u| \ge 1, {}^2\xi_u = \varnothing)\alpha_A(\lambda_1) - P(1 \le |{}^1\xi_u| \le L) - \delta_L. \end{aligned} \tag{5.9}$$

The last two terms on the right-hand side above vanish as first $u \to \infty$ and then $L \to \infty$, and, therefore, $\liminf_{t\to\infty} P({}^1\xi_t \cap A \neq \varnothing, {}^2\xi_t = \varnothing) \ge \alpha_\eta^1 \alpha_A(\lambda_1)$. Together with (5.5) this completes the proof of (5.3).

Turning to the $i = 2$ case of (5.2) and using independence of disjoint space–time regions,

$$\begin{aligned} P({}^2\xi_{t+u} \cap A \neq \varnothing) &\le P({}^2\xi_u \neq \varnothing, D_t^{(x,t+u)} \neq \varnothing \text{ for some } x \in A) \\ &\le P({}^2\xi_u \neq \varnothing) P(\zeta_t^{2,A} \neq \varnothing). \end{aligned}$$

Now let $t, u \to \infty$ to obtain

$$\limsup_{t\to\infty} P({}^2\xi_t \cap A \neq \varnothing) \le \alpha_\eta^2 \alpha_A(\lambda_2). \tag{5.10}$$

For the required lower bound we make use of the 1-blocked ancestor process $A_t^x$ defined in (4.6). By duality and independence of disjoint space–time regions, for any $L > 0$,

$$\begin{aligned} &P({}^2\xi_{t+u} \cap A \neq \varnothing) \\ &\qquad \ge \int_{|{}^2\eta|>L} P(\xi_u \in d\eta) P(\Psi(x, \bar{\xi}_t^x, \eta) = 2 \text{ for some } x \in A). \end{aligned} \tag{5.11}$$

By (4.10),

$$\begin{aligned} &P(\Psi(x, \bar{\xi}_t^x, \eta) = 2 \text{ for some } x \in A) \\ &\qquad \ge P(\bar{\xi}_t^x \neq \varnothing, A_t^x \cap {}^2\eta \neq \varnothing \text{ for some } x \in A) \\ &\qquad \ge P(\bar{\xi}_t^x \neq \varnothing \text{ for some } x \in A) - \sum_{x\in A} P(\bar{\xi}_t^x \neq \varnothing, A_t^x \cap {}^2\eta = \varnothing) \\ &\qquad = P(\zeta_t^{2,A} \neq \varnothing) - \sum_{x\in A} P(\bar{\xi}_t^x \neq \varnothing, A_t^x \cap {}^2\eta = \varnothing). \end{aligned}$$

With this bound we can appeal to (4.20), which implies that, for any $k_0 > 0$,

$$\begin{aligned} &\liminf_{t\to\infty} P(\Psi(x, \bar{\xi}_t^x, \eta) = 2 \text{ for some } x \in A) \\ &\qquad \ge \alpha_A(\lambda_2) - |A|((1-\varepsilon_1)^{k_0} + (k_0+1)\delta_L). \end{aligned}$$

Now we apply Fatou to (5.11) and obtain

$$\liminf_{t\to\infty} P({}^2\xi_{t+u}\cap A\neq\varnothing)\geq[P({}^2\xi_u\neq\varnothing)-P(1\leq|{}^2\xi_u|\leq L)]$$
$$\times[\alpha_A(\lambda_2)-|A|((1-\varepsilon_1)^{k_0}+(k_0+1)\delta_L)].$$

Finally, we let $u, L, k_0\to\infty$ in order and make use of (4.1) and (4.5) to obtain

$$\liminf_{t\to\infty} P({}^2\xi_t\cap A\neq\varnothing)\geq\alpha_\eta^2\alpha_A(\lambda_2), \tag{5.12}$$

which together with (5.10) completes the proof of (5.2) for $i=2$.

**6. Proof of Theorem 4.** Recall the notation and definitions of Theorem 4. By (5.10), (5.12) and Theorem 2, it suffices to assume that $|{}^2\eta|<\infty$ and prove

$$\begin{aligned} &P({}^1\xi_{t+u}\cap A\neq\varnothing, {}^2\xi_{t+u}=\varnothing)\\ &\qquad\to(1-\alpha_\eta^2)\mu(\zeta:\zeta\cap A\neq\varnothing)\qquad\text{as } t,u\to\infty. \end{aligned} \tag{6.1}$$

The basic idea of the proof is that when the 2's die out, they "disturb" only a bounded region of space–time. In this case, $\xi_{t+u}(x)=1$ is essentially the same as $D_{t+u}^{1,(x,t+u)}\cap{}^2\eta\neq\varnothing$, since $\lambda_1\leq\lambda^*$ implies the 1-dual is unlikely to enter the disturbed space–time region.

Let $\Gamma\subset S$ be a finite set containing ${}^2\eta$ and let $\xi_s^\Gamma$ be the multitype contact process restricted to $\Gamma$. That is, put $\xi_0^\Gamma(x)=\eta(x)$ for $x\in\Gamma$, $\xi_s^\Gamma(x)=0$ for all $s\geq 0$ and $x\in\Gamma^c$, and let the dynamics of $\xi_s^\Gamma(x)$ for $x\in\Gamma$ be the same as for $\xi_s$ except that only the Poisson processes $T^x, B^{x,y}$, $x,y\in\Gamma$ are used. Consider the events

$$\begin{aligned} E_1&=\{D_{t+u}^{1,(x,t+u)}\cap{}^1\eta\neq\varnothing\text{ for some }x\in A\},\\ E_2&=\{D_s^{1,(x,t+u)}\cap\Gamma=\varnothing\ \forall x\in A, s\in[t,t+u]\},\\ E_3&=\{{}^2\xi_u^\Gamma=\varnothing\},\\ E_4&=\{{}^2\xi_s^\Gamma={}^2\xi_s\ \forall s\in[0,u]\}. \end{aligned}$$

By independence of disjoint space–time regions, the events $E_1\cap E_2$ and $E_3$ are independent, and it is simple to check that

$$\bigcap_{i=1}^4 E_i\subset\{{}^1\xi_{t+u}\cap A\neq\varnothing, {}^2\xi_u=\varnothing\}\subset\left(\bigcap_{i=1}^4 E_i\right)\cup E_2^c\cup E_4^c. \tag{6.2}$$

We will prove (6.1) by finding appropriate estimates on the $P(E_i)$ and plugging them into (6.2).

First, by duality,

$$P(E_1) = P(\zeta^{1,A}_{t+u} \cap {}^1\eta \neq \varnothing) = P(\zeta^{1,{}^1\eta}_{t+u} \cap A \neq \varnothing) = P(\zeta_{t+u} \cap A \neq \varnothing). \tag{6.3}$$

Next, since $\{{}^2\xi_u = \varnothing\} \subset E_3 \cup E_4^c$ and $E_3 \subset \{{}^2\xi_u = \varnothing\} \cup E_4^c$,

$$P({}^2\xi_u = \varnothing) - P(E_4^c) \leq P(E_3) \leq P({}^2\xi_u = \varnothing) + P(E_4^c). \tag{6.4}$$

For fixed finite $\Gamma$, switching to forward time, $\lambda_1 \leq \lambda^*$ implies that

$$P(E_2^c) \leq \sum_{x \in A} P(\zeta^{1,x}_s \cap \Gamma \neq \varnothing \text{ for some } s \geq t) \to 0 \qquad \text{as } t \to \infty, \tag{6.5}$$

and the fact that ${}^2\eta$ is finite implies that

$$P(E_4^c) \to 0 \qquad \text{as } \Gamma \uparrow S \text{ for fixed } u > 0. \tag{6.6}$$

For a lower bound on the left-hand side of (6.1), use (6.2), (6.3) and (6.4) to obtain

$$\begin{aligned} &P({}^1\xi_{t+u} \cap A \neq \varnothing, {}^2\xi_{t+u} = \varnothing) \\ &\quad \geq P\left(\bigcap_{i=1}^3 E_i\right) - P(E_4^c) \\ &\quad = P(E_1 \cap E_2)P(E_3) - P(E_4^c) \\ &\quad \geq P(E_1)P(E_3) - P(E_2^c) - P(E_4^c) \\ &\quad \geq P(\zeta_{t+u} \cap A \neq \varnothing)P({}^2\xi_u = \varnothing) - P(E_2^c) - 2P(E_c^4). \end{aligned}$$

Since $P(\zeta_{t+u} \cap A \neq \varnothing) \to \mu(\zeta : \zeta \cap A \neq \varnothing)$ as $t \to \infty$, we let $t \to \infty$, $\Gamma \uparrow S$, and $u \to \infty$ in order above and employ (6.5) and (6.6) to obtain

$$\liminf_{t \to \infty} P({}^2\xi_t = \varnothing, {}^1\xi_t \cap A \neq \varnothing) \geq (1 - \alpha_\eta^2)\mu(\zeta : \zeta \cap A \neq \varnothing).$$

By a similar argument,

$$\begin{aligned} &P({}^1\xi_{t+u} \cap A \neq \varnothing, {}^2\xi_u = \varnothing) \\ &\quad \leq P(E_1 \cap E_2)P(E_3) + P(E_2^c) + P(E_4^c) \\ &\quad \leq P(\zeta_{t+u} \cap A \neq \varnothing)P({}^2\xi_u = \varnothing) + P(E_2^c) + 2P(E_c^4) \end{aligned}$$

and, therefore,

$$\lim_{u \to \infty} \limsup_{t \to \infty} P({}^2\xi_{t+u} \cap A \neq \varnothing, {}^1\xi_u = \varnothing) \leq (1 - \alpha_\eta^2)\mu(\zeta : \zeta \cap A \neq \varnothing).$$

Combining this with (5.4),

$$\begin{aligned} &P({}^2\xi_u \neq \varnothing, {}^2\xi_{t+u} = \varnothing) \\ &\quad \leq P({}^2\xi_u \neq \varnothing, {}^2\xi_s = \varnothing \text{ for some } s \geq u) \to 0 \qquad \text{as } u \to \infty \end{aligned}$$

gives us

$$\limsup_{t \to \infty} P({}^2\xi_t = \varnothing, {}^1\xi_t \cap A \neq \varnothing) \leq (1 - \alpha_\eta^2)\mu(\zeta : \zeta \cap A \neq \varnothing),$$

and we are done.

**Acknowledgments.** We thank Rick Durrett for showing us the proof of (1.7) and Robin Pemantle for the argument for (4.1).

Department of Mathematics
Syracuse University
Syracuse, New York 13244
USA
E-mail: jtcox@syr.edu

Department of Mathematics
University of Colorado, Colorado Springs
Colorado Springs, Colorado 80933-7150
USA
E-mail: rschinaz@uccs.edu